\documentclass{jT}

\usepackage{amssymb,amsmath,latexsym,amsthm}

\newtheorem{theorem}{Theorem}[section]
\newtheorem{lemma}[theorem]{Lemma}

\newtheorem{proposition}[theorem]{Proposition}
\newtheorem{corollary}[theorem]{Corollary}

\theoremstyle{definition}

\newtheorem*{remark}{Remark}

\begin{document}

\title[Rational approximations \ldots]{Rational approximations to $\sqrt[3]{2}$ and other algebraic numbers revisited}

\author[Paul M {\sc Voutier}]{{\sc Paul M} VOUTIER}
\address{Paul M {\sc Voutier}\\
London, UK}
\email{paul.voutier@gmail.com}

\maketitle

\begin{resume}
Dans cet article, nous am\'{e}liorons des mesures effectives d'irrationalit\'{e} 
pour certains nombres de la forme $\sqrt[3]{n}$ en utilisant des approximations obtenues
\`{a} partir de fonctions hyperg\'{e}om\'{e}triques. Ces r\'{e}sultats sont tr\`{e}s proche du mieux que peut donner cette m\'ethode. Nous obtenons ces r\'{e}sultats gr\^ace \`a des  informations arithm\'{e}tiques tr\`{e}s pr\'{e}cises sur les d\'{e}nominateurs des coefficients de ces fonctions hyperg\'{e}om\'etriques.

Des am\'eliorations  de bornes pour $\theta(k,l;x)$ et $\psi(k,l;x)$  ($k=1,3,4,6$) sont  aussi pr\'{e}sent\'{e}s.
\end{resume}

\begin{abstr}
In this paper, we establish improved effective irrationality measures
for certain numbers of the form $\sqrt[3]{n}$, using approximations obtained from
hypergeometric functions. These results are very close to the best possible using this
method. We are able to obtain these results by determining very precise arithmetic
information about the denominators of the coefficients of these hypergeometric functions.

Improved bounds for $\theta(k,l;x)$ and $\psi(k,l;x)$ for $k=1,3,4,6$ are also presented.
\end{abstr}

\bigskip
\section{Introduction}

In this article, we shall consider some refinements of a method due to Alan Baker \cite{Baker1,Baker2}
for obtaining effective irrationality measures for certain algebraic numbers of the form
$z^{m/n}$. As an example, he showed that for any integers $p$ and $q$, with $q \neq 0$, 
\begin{displaymath}
\left| 2^{1/3} - \frac{p}{q} \right| > \frac{10^{-6}}{|q|^{2.955}}. 
\end{displaymath}

This method has its basis in the work of Thue. There are two infinite 
families of hypergeometric polynomials in $\mathbb{Q}[z]$, 
${ \left\{ X_{m,n,r}(z) \right\} }_{r=0}^{\infty}$ and 
${ \left\{ Y_{m,n,r}(z) \right\} }_{r=0}^{\infty}$, such 
that ${ \left\{ Y_{m,n,r}(z)/X_{m,n,r}(z) \right\} }_{r=0}^{\infty}$ 
is a sequence of good approximations to $z^{m/n}$. 
Under certain conditions on $z$, these approximations are 
good enough to enable us to establish an effective irrationality 
measure for $z^{m/n}$ which is better than the Liouville measure. 

Since it is easy to obtain sharp estimates for the other 
quantities involved, the most important consideration in 
applying this method is the size of the denominators of 
these hypergeometric polynomials.  

Chudnovsky \cite{Chud} improved on Baker's results by showing 
that, if $p$ is a sufficiently large prime divisor of the least 
common denominator of $X_{m,n,r}(z)$ and $Y_{m,n,r}(z)$, then 
$p$ must lie in certain congruence classes mod $n$ and certain 
subintervals of $[1,nr]$. 

In the case of $z^{m/n}=2^{1/3}$, he was able to show that 
for any $\epsilon > 0$ there exists a positive integer 
$q_{0}(\epsilon)$ such that 
\begin{displaymath}
\left| 2^{1/3} - \frac{p}{q} \right| > \frac{1}{|q|^{2.4297\ldots+\epsilon}} 
\end{displaymath}
for all integers $p$ and $q$ with $|q| > q_{0}(\epsilon)$. 
Moreover, since his estimates for the relevant quantities 
are asymptotically correct, this exponent is the best that 
one can obtain from this hypergeometric method although ``off-diagonal"
or the method of ``ameliorating factors" (\`{a} la Hata) still might yield
improvements.

Shortly after this work, Easton \cite{Easton} obtained explicit versions for
the cube roots of various positive integers. For $2^{1/3}$, he showed that
\begin{displaymath}
\left| 2^{1/3} - \frac{p}{q} \right| 
> \frac{2.2 \cdot 10^{-8}}{|q|^{2.795}}
\end{displaymath}
for all integers $p$ and $q$ with $q \neq 0$.

It is the purpose of this paper to establish effective irrationality measures which
come quite close to Chudnovsky's. In the particular case of $2^{1/3}$,
\begin{displaymath}
\left| 2^{1/3} - \frac{p}{q} \right| 
> \frac{0.25}{|q|^{2.4325}} 
\end{displaymath}
for all integers $p$ and $q$ with $q \neq 0$.

This paper was initially written and circulated in 1996. Independently, Bennett \cite{Benn}
obtained a result, which in the cubic case, is slightly weaker than the theorem stated here.
E.g., for $2^{1/3}$, he showed that
\begin{displaymath}
\left| 2^{1/3} - \frac{p}{q} \right| 
> \frac{0.25}{|q|^{2.45}} 
\end{displaymath}
for all integers $p$ and $q$ with $q \neq 0$.

In fact, this subject has been the topic of even more work. 
As part of his Ph.D. Thesis (see \cite{Hei}), Heimonen has 
also obtained effective irrationality measures for numbers 
of the form $\sqrt[n]{a/b}$, as well as of the form 
$\log (a/b)$. His results are not as sharp as those of the 
author, but they are still substantially better than Easton's. 

The general method used in each of these three papers is essentially 
the same. However, there are substantial differences in the 
presentations due to the fact that the approach of Bennett and 
Heimonen shows more apparently the role that Pad\'{e} approximations 
play in this area, while the author deals explicitly with 
hypergeometric polynomials.

Actually, the referee has pointed out that other work in this area has been done,
producing results not much weaker than our own. And this work preceeded
the results of Bennett, Heimonen and the author. We are referring to the work
of Nikishin \cite{Nik} and, especially, Korobov \cite{Kor}. In particular, in 1990,
Korobov showed that
\begin{displaymath}
\left| \sqrt[3]{2} - p/q \right| > q^{-2.5},
\end{displaymath}
for all natural numbers $p$ and $q$ with $q \neq 1, 4$. The reader looking
for a more accessible reference to these works is referred to \cite[pp. 38--39]{FN}.

The main differences between this version of the paper and the previous version are
Theorem~\ref{thm:theta} and improvements in computer hardware. This has resulted in
replacing 0.93 with 0.911 in the exponents on $e$ in the expressions for $E$ and $Q$
in Theorem~\ref{thm:main} (which requires a larger value of $c_{1}$), along with the
consequent improvements to Corollary~\ref{cor:main} including new results for
$\sqrt[3]{41}$ and $\sqrt[3]{57}$.

The main incentive for publication of this paper now is completeness. Several articles
have since appeared in the literature (e.g., \cite{LPV} and \cite{TVW}) which depend on
results in this article. Furthermore, the lemmas in this article, which are either new or
sharpen results currently in the literature, are important in forthcoming articles by
the author and others. They are accompanied by an analysis showing that they are best-possible
or else what the best-possible results should be. And lastly, the main theorem itself,
along with its corollary, is an improvement on the present results in the literature.

We structure this paper as follows. Section~2 contains the statements of our results.
In Section~3, we state and prove the arithmetic results that we obtain for the 
coefficients of the hypergeometric polynomials. Section~4 is devoted to the proof
of Theorem~\ref{thm:theta}, as this theorem will be required in Section~5, where we
obtain the analytic bounds that we will require for the proof of Theorem~\ref{thm:main}.
Section~6 contains the diophantine lemma that allows us to obtain an effective
irrationality measure from a sequence of good approximations. At this point, we have
all the pieces that we need to prove Theorem~\ref{thm:main}, which is done in
Section~7. Finally, Corollary~\ref{cor:main} is proven in Section~8.

Finally, I'd like to thank Gary Walsh for his encouragement and motivation to resume
my work in this area. Also, Clemens Heuberger deserves my thanks for his careful
reading of an earlier version of this paper and accompanying suggestions. And, of
course, I thank the referee for their time and effort as well as their suggestions for
improvements.

\section{Results}

\begin{theorem}
\label{thm:main}
Let $a$ and $b$ be integers satisfying $0 < b < a$. 
Define $c_{1},d,E$ and $\kappa$ by 
\begin{eqnarray*}
d      & = & \left\{ \begin{array}{ll}
		    0   & \mbox{ if } 3 \not| \, (a-b), \\
		    1   & \mbox{ if } 3 \parallel (a-b) \mbox{ and } \\
		    3/2 & \mbox{ otherwise,}   
		\end{array}
	     \right. \\ 
E      & = & e^{-0.911} 3^{d} { \left( a^{1/2}-b^{1/2} \right) }^{-2}, \\ 
\kappa & = & \frac{\log \left\{ e^{0.911} 3^{-d} 
				{ \left( a^{1/2} + b^{1/2} \right) }^{2} 
			\right\} }
		  {\log E}
	     \mbox{ and } \\        
c_{1}  & = & 10^{40(\kappa+1)} a.
\end{eqnarray*}

If $E > 1$ then 
\begin{equation}
\label{eq:result}
\left| (a/b)^{1/3} - p/q \right| > \frac{1}{c_{1} |q|^{\kappa+1}} 
\end{equation}
for all integers $p$ and $q$ with $q \neq 0$. 
\end{theorem}

\begin{remark}
$c_{1}$ grows quite rapidly as the absolute values of the arguments of
the exponential functions in the definition of $E$ approach their best possible value
of $\pi \sqrt{3}/6=0.9068\ldots$.

In the earlier version of this paper with $0.911$ replaced by $0.93$, we could have
taken $c_{1}=10^{7(\kappa+1)} a$. It is feasible to prove Theorem~\ref{thm:main} with
0.911 replaced by 0.91, but then we would have to take $c_{1}=10^{86(\kappa+1)} a$.

The rate of growth is even more rapid as we continue to approach $0.9068$. For example,
with $0.907$, $c_{1} > 10^{2400(\kappa+1)} a$.
\end{remark}

As an application of Theorem~\ref{thm:main}, we give effective 
irrationality measures for all numbers of the form 
$\sqrt[3]{n}$ where $n$ is a cube-free rational integer  
with $2 \leq n \leq 100$ and for which the hypergeometric 
method yields an improvement over the Liouville bound. 

\begin{corollary}
\label{cor:main}
For the values of $n$ given in Table~$1$, we have 
\begin{displaymath}
\left| \sqrt[3]{n} - p/q \right| 
> \frac{c_{2}}{|q|^{\kappa+1}},  
\end{displaymath}
for all integers $p$ and $q$ with $q \neq 0$ where $c_{2}$ 
and $\kappa$ are the values corresponding to $n$ in Table~$1$. 
\end{corollary}

\begin{table}[ht]
\begin{center}
\begin{tabular}{||ccc|ccc|ccc||} \hline 
$n$   &  $c_{2}$  &  $\kappa$  & $n$   &  $c_{2}$  &  $\kappa$  & $n$   &  $c_{2}$  &  $\kappa$  \\ \hline
$2$   &  $0.25$   &  $1.4325$  & $25$  &  $0.07$   &  $1.7567$  & $60$  &  $0.08$   &  $1.5670$  \\ \hline
$3$   &  $0.37$   &  $1.6974$  & $26$  &  $0.03$   &  $1.4860$  & $61$  &  $0.06$   &  $1.5193$  \\ \hline
$4$   &  $0.41$   &  $1.4325$  & $28$  &  $0.03$   &  $1.4813$  & $62$  &  $0.04$   &  $1.4646$  \\ \hline
$5$   &  $0.29$   &  $1.7567$  & $30$  &  $0.10$   &  $1.6689$  & $63$  &  $0.02$   &  $1.3943$  \\ \hline
$6$   &  $0.01$   &  $1.3216$  & $31$  &  $0.14$   &  $1.9288$  & $65$  &  $0.02$   &  $1.3929$  \\ \hline
$7$   &  $0.08$   &  $1.6717$  & $36$  &  $0.08$   &  $1.3216$  & $66$  &  $0.04$   &  $1.4610$  \\ \hline
$9$   &  $0.08$   &  $1.6974$  & $37$  &  $0.01$   &  $1.2472$  & $67$  &  $0.06$   &  $1.5125$  \\ \hline
$10$  &  $0.15$   &  $1.4157$  & $39$  &  $0.08$   &  $1.1848$  & $68$  &  $0.08$   &  $1.5562$  \\ \hline
$11$  &  $0.22$   &  $1.8725$  & $41$  &  $0.41$   &  $1.9956$  & $70$  &  $0.12$   &  $1.6314$  \\ \hline
$12$  &  $0.28$   &  $1.9099$  & $42$  &  $0.12$   &  $1.4186$  & $76$  &  $0.08$   &  $1.5154$  \\ \hline
$13$  &  $0.35$   &  $1.8266$  & $43$  &  $0.01$   &  $1.2890$  & $78$  &  $0.03$   &  $1.5729$  \\ \hline
$15$  &  $0.19$   &  $1.4964$  & $44$  &  $0.21$   &  $1.8164$  & $83$  &  $0.09$   &  $1.6898$  \\ \hline
$17$  &  $0.01$   &  $1.1996$  & $49$  &  $0.13$   &  $1.6717$  & $84$  &  $0.37$   &  $1.8797$  \\ \hline
$18$  &  $0.37$   &  $1.9099$  & $50$  &  $0.11$   &  $1.1962$  & $90$  &  $0.09$   &  $1.3751$  \\ \hline
$19$  &  $0.02$   &  $1.2718$  & $52$  &  $0.26$   &  $1.8901$  & $91$  &  $0.009$  &  $1.2583$  \\ \hline
$20$  &  $0.009$  &  $1.1961$  & $57$  &  $0.15$   &  $1.9825$  & $98$  &  $0.38$   &  $1.4813$  \\ \hline
$22$  &  $0.07$   &  $1.2764$  & $58$  &  $0.12$   &  $1.6526$  & $100$ &  $0.35$   &  $1.4158$  \\ \hline
\end{tabular}                                                             
\caption{Results for $\sqrt[3]{n}$}
\end{center}
\end{table}

\begin{remark}
If $\alpha$ be an irrational element of $\mathbb{Q} \left( \sqrt[3]{n} \right)$, then we can write 
\begin{displaymath}
\alpha = \frac{a_{1} \sqrt[3]{n} + a_{2}}
	      {a_{3} \sqrt[3]{n} + a_{4}}, 
\end{displaymath}
where $a_{1},a_{2},a_{3},a_{4} \in \mathbb{Z}$ with 
$a_{1}a_{4}-a_{2}a_{3} \neq 0$. In this way, we can use
Corollary~\ref{cor:main} to obtain effective irrationality 
measures for any such $\alpha$ (see Section~8 of \cite{Chud}). 
\end{remark}

These values of $a$ and $b$ were found from the convergents $p/q$ in the
continued-fraction expansion of $\sqrt[3]{n}$ by setting $a/b$ to be
either $(p/q)^{3}/n$ or its reciprocal, whichever is greater than one.
For each cube-free positive integer less than or equal to $100$, we 
searched through all the convergents with $q < 10^{100}$.

In this way, we obtain measures for $\sqrt[3]{5}$, $\sqrt[3]{11}$ and $\sqrt[3]{41}$
--- values of $n$ within the range considered by Chudnovsky, but not treated by him
--- as well as an improved irrationality measure for $\sqrt[3]{7}$. Bennett also found
the same $a$ and $b$ for these $n$ (along with $n=41$ and $57$, which we also
consider here). However, his version of our Theorem~\ref{thm:main} was not
sufficiently strong to allow him to obtain effective irrationality measures for $n=41$
and $57$ which improve on Liouville's theorem, so these remain as new results here.

Given the scale of the search, the table is almost certainly complete for $n \leq 100$.

The values of $a$ and $b$ listed in Table~1 produced the minimal values
of $\kappa < 2$ satisfying the conditions of Theorem~\ref{thm:main} for
the given value of $n$.

A key element in translating the sharp result contained in Proposition~\ref{prop:denom}
into tight numerical results is a strong bound for
\begin{displaymath}
\theta(x;k,l)=\sum_{\stackrel{p \equiv l \bmod k; p, {\rm prime}}
			     {p \leq x}} 
			     \log p.
\end{displaymath}                             

Ramar\'{e} and Rumely \cite{RR} provide good bounds. However, due to recent computational work of
Rubinstein \cite{R}, we are able to improve these bounds considerably for some $k$. So we present
here the following results on $\theta(x;k,l)$, and the closely-related $\psi(x;k,l)$, for $k=1,3,4$ and $6$.

\begin{theorem}
\label{thm:theta}
\noindent
{\rm (a)}
For $1 \leq x \leq 10^{12}$,\\
$\displaystyle \max_{1 \leq y \leq x} \max \left( \left| \theta(y) - y \right|, \left| \psi(y) - y \right| \right)
\leq 2.052818 \sqrt{x}$,\\
$\displaystyle \max_{1 \leq y \leq x} \max \left( \left| \theta(y; 3, \pm 1) - y/2 \right|, \left| \psi(y; 3, \pm 1) - y \right|  \right)
\leq 1.798158 \sqrt{x}$,\\
$\displaystyle \max_{1 \leq y \leq x} \max \left( \left| \theta(y; 4, \pm 1) - y/2 \right|, \left| \psi(y; 4, \pm 1) - y \right|  \right)
\leq 1.780719 \sqrt{x}$ and \\
$\displaystyle \max_{1 \leq y \leq x} \max \left( \left| \theta(y; 6, \pm 1) - y/2 \right|, \left| \psi(y; 6, \pm 1) - y \right|  \right)
\leq 1.798158 \sqrt{x}$.

\noindent
{\rm (b)}
For each $(k,l), x_{0}$ and $\epsilon$ given in Table~$2$,
\begin{displaymath}
\left| \theta(x; k, l) - \frac{x}{\varphi(k)} \right|, \left| \psi(x; k, l) - \frac{x}{\varphi(k)} \right| \leq \epsilon x,
\end{displaymath}
for $x \geq x_{0}$.
\end{theorem}

\begin{table}[ht]
\begin{center}
\begin{tabular}{||c|c|c|c|c|c|c||} \hline 
        &  $10^{5}$  &  $10^{6}$  &  $10^{7}$  &  $10^{8}$   &   $10^{9}$  &  $10^{10}$  \\ \hline
$(1,0)$ &  $0.00474$ &  $0.00168$ & $0.000525$ & $0.0001491$ & $0.0000459$ & $0.0000186$ \\ \hline
$(3,1)$ &  $0.00405$ &  $0.00148$ & $0.000401$ & $0.0001260$ & $0.0000371$ & $0.0000351$ \\ \hline
$(3,2)$ &  $0.00217$ &  $0.00068$ & $0.000180$ & $0.0000428$ & $0.0000351$ & $0.0000351$ \\ \hline
$(4,1)$ &  $0.00494$ &  $0.00169$ & $0.000471$ & $0.0001268$ & $0.0000511$ & $0.0000511$ \\ \hline
$(4,3)$ &  $0.00150$ &  $0.00036$ & $0.000197$ & $0.0000511$ & $0.0000511$ & $0.0000511$ \\ \hline
$(6,1)$ &  $0.00405$ &  $0.00148$ & $0.000401$ & $0.0001260$ & $0.0000371$ & $0.0000351$ \\ \hline
$(6,5)$ &  $0.00217$ &  $0.00068$ & $0.000180$ & $0.0000428$ & $0.0000351$ & $0.0000351$ \\ \hline
\end{tabular}
\caption{Analytic epsilons for $x \geq x_{0}$}
\end{center}
\end{table}

Only the results for $\theta(x; 3,2)$ and $\psi(x; 3,2)$ will be used here, but we record the additional
inequalities in this theorem for use in ongoing work and by other researchers as they improve the
current bounds of Ramar\'{e} and Rumely \cite{RR} by a factor of approximately 30.

Unless otherwise noted, all the calculations mentioned in this paper were done using the
Java programming language (release 1.4.2) running on an IBM-compatible computer with
an Intel P4 CPU running at 1.8~GHz with 256 MB of memory. Source code for all programs
can be provided upon request. Many of these computations were also checked by hand,
using MAPLE, PARI/GP and UBASIC. No discrepancies beyond round-off error were found.

\section{Arithmetic Properties of Hypergeometric Polynomials}

We use ${} _{2}F_{1}(a,b;c;z)$ to denote the hypergeometric 
function
\begin{displaymath}
{} _{2}F_{1}(a,b;c;z) 
= 1 + \sum_{k=1}^{\infty} \frac{a(a+1) \cdots (a+k-1)b(b+1) \cdots (b+k-1)} 
			       {c(c+1) \cdots (c+k-1) k!} z^{k}. 
\end{displaymath}

For our purposes here, we are interested in the following functions,
which we define for all positive integers $m,n$ and $r$ with $(m,n)=1$.
Let 
\begin{eqnarray*}
X_{m,n,r}(z) & = & z^{r} {} _{2}F_{1} \left( -r, -r-m/n; 1-m/n; z^{-1} \right), \\ 
Y_{m,n,r}(z) & = & {} _{2}F_{1} \left( -r, -r-m/n; 1-m/n; z \right) \mbox{ and } \\
R_{m,n,r}(z) & = & \frac{(m/n) \cdots (r+m/n)}{(r+1) \cdots (2r+1)} 
		   {} _{2}F_{1} \left( r+1-m/n, r+1; 2r+2; 1-z \right). 
\end{eqnarray*}               

This differs from (4.3) of \cite{Chud} where the expressions 
for $X_{r}(z)$ and $Y_{r}(z)$ have been switched. The same change 
must be made in (4.4) of \cite{Chud} too. 

\vspace{3.0mm}

\noindent
{\bf Notations.} We let $D_{m,n,r}$ denote the smallest positive 
integer such that $D_{m,n,r} Y_{m,n,r}(z)$ has rational 
integer coefficients.

\noindent
To simplify the notation in the case of $m=1$ and $n=3$, which is of particular
interest in this paper, we let $X_{r}(z),Y_{r}(z),R_{r}(z)$ and $D_{r}$ denote
$X_{1,3,r}(z), Y_{1,3,r}(z), R_{1,3,r}(z)$ and $D_{1,3,r}$, respectively. 

\noindent
We will use $v_{p}(r)$ to denote the largest 
power of a prime $p$ which divides into the rational number $r$. 

\noindent
Finally, we let $\lfloor \cdot \rfloor$ denote the floor function which maps a real 
number to the greatest integer less than that number. 

\vspace{3.0mm}

We first need a refined version of Chudnovsky's Lemma~4.5 in order to establish
our criterion for the prime divisors of $D_{m,n,r}$.

\begin{lemma}
\label{lem:vp-prod}
Suppose that $m, n, p, u$ and $v$ are integers with $0<m<n$ and $(m,n)=(p,n)=1$.
For each positive integer, i, define the integer $1 \leq k_{i} \leq p^{i}$
by $k_{i}n \equiv m \bmod p^{i}$. Then
\begin{eqnarray*}
v_{p} \left( \prod_{j=u}^{v} (nj-m) \right)
& = & \sum_{i=1}^{\infty} 
      \left( \left\lfloor \frac{v-k_{i}}{p^{i}} \right\rfloor
	     - \left\lfloor \frac{u-1-k_{i}}{p^{i}} \right\rfloor \right) \\ 
& = & \sum_{i=1}^{\infty} 
      \left( \left\lfloor \frac{-u+k_{i}}{p^{i}} \right\rfloor
	     - \left\lfloor \frac{-v-1+k_{i}}{p^{i}} \right\rfloor \right).
\end{eqnarray*}
\end{lemma}

\begin{remark}
It would be more typical to state the above lemma with the condition
$0 \leq k_{i} < p^{i}$ rather than $1 \leq k_{i} \leq p^{i}$. The proof below
holds with either condition. However, the above formulation suits our needs
in  the proof of Proposition~\ref{prop:denom} below better.
\end{remark}

\begin{proof}
For each positive integer $i$, we will count the number of $j$'s in
$u \leq j \leq v$ with $nj-m \equiv 0 \bmod p^{i}$. That is, with
$nj-k_{i}n \equiv 0 \bmod p^{i}$. And, since $(n,p)=1$, with
$j \equiv k_{i} \bmod p^{i}$. The remainder of the proof is identical
to Chudnovsky's proof of his Lemma~4.5 \cite{Chud}, upon replacing
his $p$ with $p^{i}$.
\end{proof}

\begin{proposition}
\label{prop:denom}
Let $m,n$ and $r$ be positive integers with $0 < m < n$ and 
$(m,n)=1$.

The largest power to which a prime $p$ can divide $D_{m,n,r}$ is at most the number
of positive integers $i$ for which there exist a positive integer $l_{i}$ satisfying
$(l_{i},n)=1, l_{i} p^{i} \equiv -m \bmod n$ such that
\begin{displaymath}
\frac{l_{i}p^{i}+m}{n} \leq r \bmod p^{i} \leq \frac{(n-l_{i})p^{i}-m-n}{n}.
\end{displaymath}
Furthermore, all such $i$ satisfy $p^{i} \leq nr$.
\end{proposition}

\begin{remark}
From the calculations done in the course of this, and other, work (see, for
example, the notes following Lemmas~\ref{lem:denom-gen}, \ref{lem:denom-3}
and \ref{lem:coeff-ub}), it appears that the conditions given in this Proposition
provide the exact power to which a prime divides $D_{m,n,r}$.
However, I have not been able to prove this. 
\end{remark}

\begin{proof}
Let $a_{r,h}$ denote the coefficient of $z^{h}$ in 
$Y_{m,n,r}(z)$ and let $p$ be a prime number. From our 
definition of $Y_{m,n,r}(z)$ above, we can write 
\begin{displaymath}
a_{r,h} = {r \choose h} \frac{C_{r,h}}{B_{r,h}},
\end{displaymath}
where 
\begin{displaymath}
B_{r,h} = \prod_{i=1}^{h} (in-m) 
\hspace{3.0mm} \mbox{ and } \hspace{3.0mm}
C_{r,h} = \prod_{i=r-h+1}^{r} (in+m).  
\end{displaymath}

We first show that if $p$ divides $D_{m,n,r}$ then $(p,n)=1$. 

If $p$ does divide $D_{m,n,r}$ then $p$ must divide $B_{r,h}$
for some $0 \leq h \leq r$. So it must divide some number of the
form $in-m$ where $1 \leq i \leq r$. But, if $p$ divides such a
number and also divides $n$, then it must also divide $m$. However, 
our hypothesis that $(m,n)=1$ does not allow this and so,  
if $p$ divides $D_{m,n,r}$ then $(p,n)=1$. 

Therefore, for any positive integer $i$, we can find an integer $k_{i}$ with
$1 \leq k_{i} \leq p^{i}, (k_{i},p^{i})=1$ and  $k_{i}n \equiv m \bmod p^{i}$.

As $1 \leq k_{i}$ and $m<n$, we know that $0 < k_{i}n-m$, and so there
must be a positive integer $l_{i}$ with $(l_{i},n)=1$ and $k_{i}n-m=l_{i}p^{i}$.
Furthermore, $l_{i} < n$.

Returning to our expression for $a_{r,h}$, we have 
\begin{displaymath}
v_{p} \left( a_{r,h} \right)
= v_{p} \left( {r \choose h} \right) + v_{p} \left( C_{r,h} \right)
- v_{p} \left( B_{r,h} \right).
\end{displaymath}

It is well-known that
\begin{displaymath}
v_{p} \left( {r \choose h} \right) 
= \sum_{i=1}^{\infty} 
  \left( \left\lfloor \frac{r}{p^{i}} \right\rfloor
         - \left\lfloor \frac{h}{p^{i}} \right\rfloor
         - \left\lfloor \frac{r-h}{p^{i}} \right\rfloor \right).
\end{displaymath}

From the first expression in Lemma~\ref{lem:vp-prod} with $u=1$ and $v=h$,
\begin{displaymath}
v_{p} \left( B_{r,h} \right) 
= \sum_{i=1}^{\infty} 
  \left( \left\lfloor \frac{h-k_{i}}{p^{i}} \right\rfloor
         - \left\lfloor \frac{-k_{i}}{p^{i}} \right\rfloor \right).
\end{displaymath}

From the second expression in Lemma~\ref{lem:vp-prod} with $u=-r$ and $v=-r+h-1$,
\begin{displaymath}
v_{p} \left( C_{r,h} \right) 
= \sum_{i=1}^{\infty} 
  \left( \left\lfloor \frac{r+k_{i}}{p^{i}} \right\rfloor
         - \left\lfloor \frac{r+k_{i}-h}{p^{i}} \right\rfloor \right).
\end{displaymath}

Thus, we want to determine when 
\begin{equation}
\label{eq:positive}
\left\lfloor \frac{r}{p^{i}} \right\rfloor - \left\lfloor \frac{h}{p^{i}} \right\rfloor 
- \left\lfloor \frac{r-h}{p^{i}} \right\rfloor
- \left\lfloor \frac{h-k_{i}}{p^{i}} \right\rfloor + \left\lfloor \frac{-k_{i}}{p^{i}} \right\rfloor 
+ \left\lfloor \frac{r+k_{i}}{p^{i}} \right\rfloor - \left\lfloor \frac{r+k_{i}-h}{p^{i}} \right\rfloor 
\end{equation}
is negative.

This will suffice for the purpose of proving this proposition since, as we shall
show shortly, the expression in (\ref{eq:positive}) can never be less than $-1$. 

We now show that if $p^{i} > nr$, then the expression in (\ref{eq:positive})
cannot be negative. This will establish the last statement in the Proposition.

Since $0 \leq h \leq r < p^{i}$ for such $i$, the first three terms in
(\ref{eq:positive}) are 0. Furthermore, the same inequalities for $h$
and $r$ along with the fact that $k_{i} > 0$ show that the sum of the
last two terms cannot be negative.

We saw above that $k_{i}n-m=l_{i}p^{i}$ for a positive integer $l_{i}$.
So it follows that $k_{i}n \geq p^{i}+m > nr \geq nh$. In particular,
$k_{i} > h$. Furthermore, $1 \leq k_{i} \leq p^{i}$. Therefore,
$\lfloor (h-k_{i})/p^{i} \rfloor$ and $\lfloor -k_{i}/p^{i} \rfloor$, are both
equal to $-1$, so the sum of the remaining terms in (\ref{eq:positive})
is also zero.

This establishes the last statement in the Proposition.

Moreover, if $p^{i}>nr+m$, then
\begin{displaymath}
r+k_{i} < \frac{p^{i}-m}{n}+ \frac{(n-1)p^{i}+m}{n} < p^{i}.
\end{displaymath}

And so the expression in (\ref{eq:positive}) is always 0 for such $i$.
We will use this fact in the proof of Lemmas~\ref{lem:denom-gen}
and \ref{lem:denom-3} below.

For any positive integer $i$, we can write $h$ and $r$ uniquely as
\begin{displaymath}
h = h_{i1}p^{i}+h_{i0} \hspace{3.0mm} \mbox{ and } \hspace{3.0mm}
r = r_{i1}p^{i}+r_{i0}, 
\end{displaymath}
where $0 \leq h_{i0},r_{i0} < p^{i}$.  

With this notation, we see that 
\begin{eqnarray}
\label{eq:others}
\left\lfloor \frac{r}{p^{i}} \right\rfloor - \left\lfloor \frac{h}{p^{i}} \right\rfloor 
- \left\lfloor \frac{r-h}{p^{i}} \right\rfloor
& = & - \left\lfloor \frac{r_{i0}-h_{i0}}{p^{i}} \right\rfloor 
      \nonumber \\
\left\lfloor \frac{h-k_{i}}{p^{i}} \right\rfloor - \left\lfloor \frac{-k_{i}}{p^{i}} \right\rfloor
& = & h_{i1} + 1 + \left\lfloor \frac{h_{i0}-k_{i}}{p^{i}} \right\rfloor
      \mbox{ and } \\
\left\lfloor \frac{r+k_{i}}{p^{i}} \right\rfloor - \left\lfloor \frac{r+k_{i}-h}{p^{i}} \right\rfloor
& = & h_{i1} + \left\lfloor \frac{r_{i0}+k_{i}}{p^{i}} \right\rfloor
	    - \left\lfloor \frac{r_{i0}+k_{i}-h_{i0}}{p^{i}} \right\rfloor.
      \nonumber
\end{eqnarray}

The first relation holds since $1 \leq h_{i0}, r_{i0} < p^{i}$ and so
$\lfloor h_{i0}/p^{i} \rfloor = \lfloor r_{i0}/p^{i} \rfloor =0$. The second
relation holding since $1 \leq k_{i} \leq p^{i}$ and so $\lfloor -k_{i}/p^{i} \rfloor =-1$. 

The last two quantities can only have the values $h_{i1}$ or $h_{i1}+1$,
so if the expression in (\ref{eq:positive}) is to be negative then the first
quantity here must be zero, since it is never negative, the second must
be $h_{i1}+1$ and the third must be $h_{i1}$. This information also
substantiates our claim above that the expression in (\ref{eq:positive})
is always at least $-1$. 

Since $0 \leq h_{i0}, r_{i0} < p^{i}$, the first quantity in
(\ref{eq:others}) is zero if and only if 
\begin{equation}
\label{eq:rh}
r_{i0} \geq h_{i0}. 
\end{equation}

The second quantity in (\ref{eq:others}) is $h_{i1}+1$ if and only if
\begin{equation}
\label{eq:hk}
h_{i0} \geq k_{i}. 
\end{equation}

Finally, if the last quantity in (\ref{eq:others}) is $h_{i1}$, then
$\lfloor (r_{i0}+k_{i})/p^{i} \rfloor = \lfloor (r_{i0}+k_{i}-h_{i0})/p^{i} \rfloor$.
From (\ref{eq:hk}), we find that $r_{i0}+k_{i}-h_{i0} \leq r_{i0} 
< p^{i}$, so $r_{i0}+k_{i} < p^{i}$ also. Hence 
\begin{equation}
\label{eq:rk}
0 < \frac{r_{i0}+k_{i}}{p^{i}} < 1, 
\end{equation}
the left-hand inequality being strict since $k_{i} > 0$. 

From (\ref{eq:rh}), we have $k_{i} \leq r_{i0}+k_{i}-h_{i0}$, 
while from (\ref{eq:hk}) and (\ref{eq:rk}), it follows that  
$r_{i0}+k_{i}-h_{i0} < p^{i}-h_{i0} \leq p^{i}-k_{i}$. Combining 
these inequalities, we find that 
\begin{displaymath}
\frac{k_{i}}{p^{i}} \leq \frac{r_{i0}+k_{i}-h_{i0}}{p^{i}} 
< 1 - \frac{k_{i}}{p^{i}}. 
\end{displaymath}

In addition, from (\ref{eq:rh}) and (\ref{eq:hk}), 
we know that $k_{i} \leq r_{i0}$ and from (\ref{eq:rk}), 
$r_{i0} \leq p^{i}-k_{i}-1$, so 
\begin{displaymath}
r_{i1} + \frac{k_{i}}{p^{i}} \leq \frac{r}{p^{i}} 
\leq r_{i1} + 1 - \frac{k_{i}+1}{p^{i}}. 
\end{displaymath}

Substituting $(l_{i}p^{i}+m)/n = k_{i}$ into this expression
completes the proof of the Proposition.
\end{proof}

It will be helpful for applications to present a slightly weaker but more
immediately applicable result on the prime divisors of $D_{m,n,r}$. With
that in mind, we state the following.

\begin{lemma}
\label{lem:denom-gen}
{\rm (a)} Let $r$ be a positive integer. If 
$p|D_{m,n,r}$, then
\begin{displaymath}
v_{p} \left( D_{m,n,r} \right) 
\leq \left\lfloor \frac{\log (nr)}{\log p} \right\rfloor.
\end{displaymath}

\noindent
{\rm (b)} If $p$ is a prime number greater than $(nr)^{1/2}$ 
which is a divisor of $D_{m,n,r}$, then $p^{2} \not| D_{m,n,r}$
and for some $1 \leq l < n/2$ with $(l,n)=1$, $lp \equiv -m \bmod n$, 
and 
\begin{equation}
\label{eq:pcond-gen}
\frac{nr+m+n}{nA+n-l} \leq p \leq \frac{nr-m}{nA+l},  
\end{equation}
for some non-negative integer $A$. Moreover, every such prime 
greater than $(nr+m)^{1/2}$ is a divisor of $D_{m,n,r}$.
\end{lemma}

\begin{remark}
The result in (a) is best possible. E.g., $D_{2,3,17}$ is divisible by $49$
and $\lfloor \log (3 \cdot 17)/(\log 7) \rfloor =2$. This example also shows
that neither of the statements in Lemma~\ref{lem:denom-gen} holds here
(i.e., smaller primes may divide $D_{m,n,r}$ to higher powers than one
and $p$ need not lie in the intervals specified by (\ref{eq:pcond-gen})).

Furthermore, $5$ divides $D_{2,3,10}$, so the congruence
conditions in Lemma~\ref{lem:denom-gen} do not hold in general either
(since $1 \cdot 5 \not\equiv -2 \bmod 3$).
\end{remark}

\begin{proof}
(a) This follows immediately from the last statement in Proposition~\ref{prop:denom}.

(b) Again from the last statement in Proposition~\ref{prop:denom} and our
lower bound for $p$, we need only consider $i=1$.

From the inequality on $r \bmod p$ in Proposition~\ref{prop:denom},
we can write
\begin{equation}
\label{eq:rmodp-gen}
Ap + \frac{lp+m}{n} \leq r \leq Ap + \frac{(n-l)p-m-n}{n},
\end{equation}
for some non-negative $A$. This provides our upper and lower bounds
for $p$ in part~(b), which suffices to prove the  first statement
in part~(b). 

To prove the second statement, we will show that these primes divide
the denominator of the leading coefficient of $Y_{r}(z)$. So we let
the quantity denoted by $h$ in the proof of Proposition~\ref{prop:denom}
be $r$. Using the arguments to derive (\ref{eq:others}) in the proof
of Proposition~\ref{prop:denom}, (\ref{eq:positive}) simplifies to
\begin{displaymath}
-1 - \left\lfloor \frac{r_{i0}-k_{i}}{p^{i}} \right\rfloor
+ \left\lfloor \frac{r_{i0}+k_{i}}{p^{i}} \right\rfloor - \left\lfloor \frac{k_{i}}{p^{i}} \right\rfloor \\
= -1 - \left\lfloor \frac{r_{i0}-k_{i}}{p^{i}} \right\rfloor + \left\lfloor \frac{r_{i0}+k_{i}}{p^{i}} \right\rfloor,
\end{displaymath}
where $r \equiv r_{i0} \bmod p^{i}$.

Therefore, as we saw in the proof of Proposition~\ref{prop:denom},
\begin{displaymath}
v_{p}(a_{r,r}) = \sum_{i=1}^{\infty} \left( -1 - \left\lfloor \frac{r_{i0}-k_{i}}{p^{i}} \right\rfloor
                                            + \left\lfloor \frac{r_{i0}+k_{i}}{p^{i}} \right\rfloor \right).
\end{displaymath}

Notice that for $i \geq 2$, $p^{i} > nr+m$, so, as we saw in the proof of
Proposition~\ref{prop:denom}, the summands for such $i$ are zero and can be ignored.

For $i=1$, from (\ref{eq:rmodp-gen}), $(lp+m)/n \leq r_{i0} \leq ((n-l)p-m-n)/n$.
From the relationship between $k_{1}$ and $l_{1}$ given in the proof of
Proposition~\ref{prop:denom}, we also have $k_{1}=(lp+m)/n$. Therefore,
$0 \leq \lfloor r_{10}-k_{1} \rfloor, \lfloor r_{10}+k_{1} \rfloor \leq p-1$ and the
summand for $i=1$ is -1. Hence $v_{p}(a_{r,r})=-1$, so $p$ divides the denominator
of $a_{r,r}$ precisely once, completing the proof of the Lemma.
\end{proof}

As $n$ gets larger, the structure of the denominator becomes more complicated and the
above is the best that we can do. However, in the case of $m=1$ and $n=3, 4$ or $6$,
we can obtain a sharper result which will be used in this paper.

\begin{lemma}
\label{lem:denom-3}
Let $m=1$ and $n=3, 4$ or $6$.

\noindent
{\rm (a)} Let $r$ be a positive integer. If 
$p|D_{m,n,r}$ then $p \equiv n-1 \bmod n$ and 
\begin{displaymath}
v_{p} \left( D_{m,n,r} \right) 
\leq \left\lfloor \frac{\log (nr)}{2 \log p} + \frac{1}{2} \right\rfloor.
\end{displaymath}

\noindent
{\rm (b)} If $p$ is a prime number greater than $(nr)^{1/3}$ 
which is a divisor of $D_{m,n,r}$ then $p \equiv n-1 \bmod n$, 
$p^{2} \not| D_{m,n,r}$ and 
\begin{equation}
\label{eq:pcond-3}
\frac{nr+n+1}{nA+n-1} \leq p \leq \frac{nr-1}{nA+1},  
\end{equation}
for some non-negative integer $A$. Moreover, every such prime 
greater than $(nr+1)^{1/2}$ is a divisor of $D_{m,n,r}$.
\end{lemma}

\begin{remark}
The result in (a) is best possible. $D_{42}$ is divisible by 25 and
$\lfloor \log (3 \cdot 42)/(2 \log 5)$  $+ 1/2 \rfloor =2$.
Similarly, $D_{1042}$ is divisible by 125 and
$\lfloor \log (3 \cdot 1042)/(2 \log 5) + 1/2 \rfloor = 3$. However, it is not
true that $v_{p}(D_{r}) \geq 2$ for all $p \leq(3r)^{1/3}$ (e.g., $v_{5}(D_{43})=1$).
\end{remark}

\begin{remark}
The second statement in (b) holds for all $p > (nr)^{1/3}$, and this
is best possible as the example in the previous remark shows,
however the proof is technical and lengthy. Furthermore, the result
here suffices for our needs below.
\end{remark}

\begin{proof}
(a) We apply Proposition~\ref{prop:denom}. As we saw there, $(p,n)=1$.
For these values of $n$, the only integers less than $n$ and relatively prime
to $n$ are $1$ and $n-1$.

If $p \equiv 1 \bmod n$ or if $p \equiv n-1 \bmod n$ and $i$ is even,
then we require $l_{i} \equiv n-1 \bmod n$ to satisfy $l_{i}p^{i} \equiv -1 
\bmod n$. However, with this value of $l_{i}$,
\begin{displaymath}
\frac{(n-1)p^{i}+1}{n} = \frac{l_{i}p^{i}+m}{n}
\leq r \bmod p^{i}
\leq \frac{(n-l_{i})p^{i}-m-n}{n} = \frac{p^{i}-n-1}{n}
\end{displaymath}
can never be satisfied.

If $p \equiv n-1 \bmod n$ and $i$ is odd, then we can take $l_{i}=1$.
From the last statement in Proposition~\ref{prop:denom}, $p^{i} \leq nr$,
so the largest possible $i$ is at most $\log(nr)/\log(p)$, a fact which
completes the proof of part~(a).

(b) The same argument as for $p\ \equiv n-1 \bmod n$ in the proof
of part~(a) shows that we need only consider $i=1$ for $p > (nr)^{1/3}$.

The remainder of the proof is identical to the proof of Lemma~\ref{lem:denom-gen}(b).
\end{proof}

\begin{lemma}
\label{lem:numer}
Let $m,n$ and $r$ be positive integers with $(m,n)=1$. 
Define $\mu_{n} = \prod_{p|n} p^{1/(p-1)}$ and $s_{n,r} 
= \prod_{p|n} p^{v_{p}(r!)}$.  

\noindent
{\rm (a)} Let $d$ be a positive divisor of $n$.
The numerators of the coefficients of the polynomials 
$X_{m,n,r}(1-dz)$ and $Y_{m,n,r}(1-dz)$ are divisible by $d^{r}$. 

\noindent
{\rm (b)} The numerators of the coefficients of the polynomials 
$X_{m,n,r} \left( 1 - n\mu_{n}z \right)$ and 
$Y_{m,n,r} \left( 1-n\mu_{n}z \right)$ are divisible by $n^{r} s_{n,r}$. 
\end{lemma}

\begin{proof}
(a) This is a variation on part~(b) which will prove useful both 
here and elsewhere. Its proof is virtually identical to the proof of part~(b).

(b) This is Proposition~5.1 of \cite{Chud}.  
\end{proof}

\section{Proof of Theorem~\ref{thm:theta}}

(a) The bounds for $x \leq 10^{12}$ are determined through direct calculation.
We coded the Sieve of Eratosthenes in Java and ran it, in segments of size $10^{8}$,
to determine all primes less than $10^{12}$ as well as upper and lower bounds for
$\theta(x; k, l)$ and $\psi(x; k, l)$ for $x \leq 10^{12}$. The entire computation
took approximately $182,000$ seconds.

As Ramar\'{e} and Rumely note, considerable roundoff error can arise in the sum of
so many floating point numbers. We handled this issue in a similar way to them.
We multiply each log by $10^{6}$, round the resulting number down to the greatest
integer less than the number as a lower bound and round it up to the least integer
greater than the number as an upper bound. We then sum these integers and store the
sums in variables of type long, which have a maximum positive value of $2^{63}-1 = 9.233... \cdot 10^{18}$
-- a number greater than our sums. This is more crude than Ramar\'{e} and Rumely's
method, but sufficiently accurate for our needs here.

In addition to just establishing the desired inequalities, we also compute, and have stored,\\
(i) our upper and lower bounds for $\theta(10^{8}i; k, l)$ and  $\psi(10^{8}i; k, l)$,\\
(ii) $\pi(10^{8}i, k, l)$,\\
(iii) $\displaystyle \min_{x \in (10^{8}(i-1), 10^{8}i]} \frac{\theta(x; k,l) - x/\varphi(k)}{\sqrt{x}},$
\hspace{3.0mm} $\displaystyle \min_{x \in (10^{8}(i-1), 10^{8}i]} \frac{\psi(x; k,l) - x/\varphi(k)}{\sqrt{x}}$ and \\
(iv) $\displaystyle \max_{x \in (10^{8}(i-1), 10^{8}i]} \frac{\theta(x; k,l) - x/\varphi(k)}{\sqrt{x}},$
\hspace{3.0mm} $\displaystyle \max_{x \in (10^{8}(i-1), 10^{8}i]} \frac{\psi(x; k,l) - x/\varphi(k)}{\sqrt{x}}$,\\
for $i=1,...,10,000$.

(b) The bounds in part (b) are obtained by applying Theorem~5.1.1 of \cite{RR} with the L-function zero
information calculated by Michael Rubinstein \cite{R}. We include details of the values
used in Table~3, where we round all quantities up by one in the seventh significant decimal
(sixth decimal for $\tilde{A}_{\chi}$, $\tilde{B}_{\chi}$, $\tilde{C}_{\chi}$, $\tilde{D}_{\chi}$
for the sake of space).

Note that for these values of $k$, there is only one character, $\chi$, for each $d$.

For the computation of $\tilde{A}_{\chi}$, we followed the advice of Ramar\'{e} and
Rumely \cite[p.~414]{RR} regarding the evaluation of their $K_{1}$ and $K_{2}$.
Using Simpson's rule with an interval size of $0.001$ (along with their Lemma~4.2.4), we bound from
above the integral for $K_{n}(z,w)$ in their equation~(4.2.4) for $u=w \ldots 1000$. We then apply
their Lemma~4.2.3 with $w=1000$, which is sufficiently large to provide a good upper bound.

This provides us with an upper bound for $\epsilon(\psi, x, k)$.

Using the authors' upper bound for $\epsilon(\theta, x, k)$ on page~420 of \cite{RR},
we see that our results holds for $x \geq x_{0}$.

Proceeding as above, we found agreement with the data that Ramar\'{e} and Rumely present in their Table~1
for $k=1, 3$ and $4$.

\begin{table}[ht]
\begin{center}
\begin{tabular}{||c|c|c|c|c|c||} \hline 
$k$                    &      $1$                  & \multicolumn{2}{|c|}{$3$}                           & \multicolumn{2}{|c||}{$4$}                           \\ \hline
$m$                    & $14$                      & \multicolumn{2}{|c|}{$14$}                          & \multicolumn{2}{|c||}{$14$}                          \\ \hline
$\delta$               & $6.289071 \cdot 10^{-7}$  & \multicolumn{2}{|c|}{$1.256642 \cdot 10^{-6}$}      & \multicolumn{2}{|c||}{$1.798450 \cdot 10^{-6}$}      \\ \hline
$A(m, \delta)$         & $1.082027 \cdot 10^{91}$  & \multicolumn{2}{|c|}{$6.691384 \cdot 10^{86}$}      & \multicolumn{2}{|c||}{$4.425147 \cdot 10^{84}$}      \\ \hline
$\tilde{R}$            & $2.721552 \cdot 10^{-11}$ & \multicolumn{2}{|c|}{$9.085095 \cdot 10^{-11}$}     & \multicolumn{2}{|c||}{$1.207835 \cdot 10^{-10}$}     \\ \hline
$\epsilon(\psi, x, k)$ & $3.613190 \cdot 10^{-5}$  & \multicolumn{2}{|c|}{$7.097148 \cdot 10^{-5}$}      & \multicolumn{2}{|c||}{$1.001340 \cdot 10^{-4}$}      \\ \hline
$d$                    & $1$                       & $1$                      & $3$                      &  $1$                      & $4$                      \\ \hline
$H_{\chi}$             & $8000000.365$             & $4000000.042$            & $4000000.413$            & $2800000.0623$            & $2800000.340$            \\ \hline
$\tilde{A}_{\chi}$     & $5.81243 \cdot 10^{-98}$  & $8.94572 \cdot 10^{-94}$ & $9.83501 \cdot 10^{-94}$ & $1.27527 \cdot 10^{-91}$ &  $1.43730 \cdot 10^{-91}$ \\ \hline
$\tilde{B}_{\chi}$     & $9.09392 \cdot 10^{-103}$ & $2.85005 \cdot 10^{-98}$ & $3.04716 \cdot 10^{-98}$ & $5.86164 \cdot 10^{-96}$ &  $6.37182 \cdot 10^{-96}$ \\ \hline
$\tilde{C}_{\chi}$     & $7.30396 \cdot 10^{-98}$  & $1.13798 \cdot 10^{-93}$ & $1.23103 \cdot 10^{-93}$ & $1.63333 \cdot 10^{-91}$ &  $1.80646 \cdot 10^{-91}$ \\ \hline
$\tilde{D}_{\chi}$     & $1.14495 \cdot 10^{-102}$ & $3.63318 \cdot 10^{-98}$ & $3.82113 \cdot 10^{-98}$ & $7.52419 \cdot 10^{-96}$ &  $8.02375 \cdot 10^{-96}$ \\ \hline
$\tilde{E}_{\chi}$     & $31.414915$               & $28.3898896$              & $33.1560902$            & $ 26.8928884$             &  $32.8584828$            \\ \hline
\end{tabular}
\caption{Data for the Proof of Theorem~\ref{thm:theta}}
\end{center}
\end{table}

Since $2.052818 < 0.0000186x^{1/2}$ for $x \geq 12.2 \cdot 10^{9}$, the stated inequalities for
$\theta(x)$ and  $\psi(x)$ holds for such $x$. Using the above sieve code, it is straightforward to calculate
$\theta(x)$ and  $\psi(x)$ for $x < 12.2 \cdot 10^{9}$. These calculations complete the proof
of (b) for $\theta(x)$ and  $\psi(x)$.

Similarly, $1.798158 < 0.0000351x^{1/2}$ for $x \geq 2.7 \cdot 10^{9}$ and a computation completes the proof of (b) for $k=3$ and 6.

Finally, $1.780719 < 0.0000511x^{1/2}$ for $x \geq 2.7 \cdot 10^{9}$ and a computation completes the proof of (b) for $k=4$..

\section{Analytic Properties of Hypergeometric Polynomials}

\begin{lemma}
\label{lem:coeff-ub}
Let $r$ be a positive integer and define $N_{r}$ to be the greatest common
divisor of the numerators of the coefficients of $X_{r}(1-(a-b)x/a)$, where
$a, b$ and $d$ are as defined in Theorem~$\ref{thm:main}$.

\noindent
{\rm (a)} We have
\begin{displaymath}
\frac{1}{200} 
  <   \frac{0.29D_{r}^{2}}{r^{1/6} 4^{r}}.
\end{displaymath}

\noindent
{\rm (b)} We have
\begin{eqnarray*}
\frac{3^{dr} D_{r}}{N_{r}} 
  <   1.161 \cdot 10^{39} e^{0.911r} 
\mbox{ and }
\frac{(1/3) \cdots (r+1/3)}{r!} \frac{3^{dr} D_{r}}{N_{r}} 
< 1.176 \cdot 10^{40} e^{0.911r}.
\end{eqnarray*}
\end{lemma}

\begin{remark}
These results are very close to best possible. 
Chudnovsky \cite{Chud} has shown that  
$D_{r} \sim e^{\pi \sqrt{3}\, r/6}= e^{0.9068\ldots r}$ 
as $r \rightarrow \infty$. 
\end{remark}

\begin{remark}
We were able to calculate the $D_{r}$ exactly for all $r \leq 2000$
(with two different methods using both Java and UBASIC~8.8). These actual values were
equal to the values calculated using Proposition~\ref{prop:denom}. This strengthens
our belief that Proposition~\ref{prop:denom} captures the precise behaviour of the
prime divisors of $D_{m,n,r}$ (at least for $m=1, n=3$).
\end{remark}

\begin{proof}
We will establish both parts of this lemma via computation 
for $r$ up to the point where Theorem~\ref{thm:theta} can be used
to prove the lemma for all larger $r$.

(a) We computed the quantity on the right-hand side for all $r \leq 2000$,
as part of the computation for part (b). We found that its minimum is
$0.00501 \ldots$, which occurs at $r=13$.

From the second statement in Lemma~\ref{lem:denom-3}(b), we know that if
$p$ is a prime congruent to $2 \bmod 3$ with $(3r+4)/2 \leq p \leq 3r-1$,
then $p | D_{r}$. Since we may now assume that $r > 2000$, we know that
$(3r+4)/2 > 3000$.

From Theorem~\ref{thm:theta} and a bit of computation, for $x > 3000$, we
find that $|\theta(x;3,2)-x/2|<0.011x$, so the product of the primes congruent
to $2 \bmod 3$ in that interval is at least $e^{0.7r-1.511}$. Therefore,
$D_{r}/4^{r}>e^{0.014r-1.511}$. Since $r^{1/6} = e^{(\log r)/6}<e^{0.0007r}$
for $r \geq 2000$, the desired result easily follows.

(b) Here the computation needs to include much larger values of $r$,
so we need to proceed more carefully.

We break the computation into several parts.

(1) The computation of the factorial and factorial-like product on the left-hand side of the
second inequality. We shall see below that the product of these terms grows quite slowly and
they have a simple form, so this computation is both easy and fast.

(2) The computation of $3^{dr}/N_{r}$. From Lemma~\ref{lem:numer}, we find that if $d=0$ 
then $\left( 3, N_{r} \right) = 1$, if $d=1$ then $3^{r} | N_{r}$ and if $d=3/2$ then
$3^{r+v_{3}(r!)} | N_{r}$. For $d=3/2$, one can often do better, by directly calculating
the numerators of the coefficients of $X_{r}(1-3\sqrt{3}\, x)$, by means of equations
(5.2)--(5.4) in the proof of Chudnovsky's Proposition~5.1 \cite{Chud}.

Directly calculating $N_{r}$ is substantially more time-consuming than calculating 
$3^{r+v_{3}(r!)}$, so we always calculate $3^{r+v_{3}(r!)}$, continue with
calculating $D_{r}$ and only perform the direct calculation of $N_{r}$ if the
size of $3^{3r/2}D_{r}/3^{r+v_{3}(r!)}$ warrants it.

(3) The computation of the contribution to $D_{r}$ from the small primes, that is those
less than $\sqrt[3]{3r}$, using Proposition~\ref{prop:denom}.

To speed up this part of the calculation, and the following parts, the primes and their
logarithms do not have to be recalculated for each $r$. Instead, we calculate and store
the first million primes congruent to 2 mod 3 (the last one being 32,441,957) and
their logarithms before we start the calculations for any of the $r$'s.

(4) The computation of the contribution to $D_{r}$ of all primes from $\sqrt[3]{3r}$
to $(3r-1)/(3A(r)+1)$ for some non-negative integer $A(r)$, which depends only on $r$.
Again, we use Proposition~\ref{prop:denom} as well as the cached primes and their
logarithms here.

(5) The computation of the contribution to $D_{r}$ from the remaining larger primes.

From Lemma~\ref{lem:denom-3}(b), we can see that for any non-negative integer $A$,
the contribution to $D_{r}$ from the primes satisfying (\ref{eq:pcond-3}) changes,
as we increment $r$, by at most the addition of the log of one prime, if there 
is a prime congruent to $2 \bmod 3$ between  $3(r-1)/(3A+1)$ and $3r/(3A+1)$,
and the subtraction of another, if there is a prime congruent to $2 \bmod 3$
between $3(r-1)/(3A+2)$ and $3r/(3A+2)$. This fact makes it very quick to compute
the contribution from these intervals for $r$ from the contribution from these intervals
for $r-1$ --- much quicker than recomputing them directly. So we incorporate this strategy
here: for each $i < A(r)$, we store the smallest and largest primes in these intervals
along with the sum of the logarithms of the primes, $p \equiv 2 \bmod 3$, in these 
intervals.

Again, we use the cached primes and their logarithms for the intervals that lie within the cache.

In this manner, we proceeded to estimate the size of the required quantities for 
all $r \leq 200,000,000$. This computation took approximately 89,000 seconds.

The maximum of $3^{dr} D_{r}/ \left( N_{r} e^{0.911r} \right)$ occurs at
$r=19,946$ and is less than $1.161 \cdot 10^{39}$, while the maximum of
$(1/3) \cdots (r+1/3) 3^{dr} D_{r} / \left( N_{r} e^{0.911r} r! \right)$
also occurs at $r=19,946$ and is less than $1.176 \cdot 10^{40}$.  

For $r > 200 \cdot 10^{6}$, we can use the analytic estimates in Theorem~\ref{thm:theta}.

From Lemma~\ref{lem:numer}, we know that $3^{dr} N_{r}^{-1} \leq 3^{r/2-v_{3}(r!)}$.
In addition, $r/2-v_{3}(r!) \leq (\log r)/(\log 3) + 0.5$ and 
\begin{displaymath}
\frac{(1/3) \cdots (r+1/3)}{r!} 
\leq \frac{4}{9} \exp \left( \int_{1}^{r} \frac{dx}{3x} \right) 
\leq \frac{4r^{1/3}}{9}, 
\end{displaymath}
for $r \geq 1$, so 
\begin{equation}
\label{eq:fact}
\frac{3^{dr}}{N_{r}} < 1.8r \hspace{3.0mm} \mbox{ and } \hspace{3.0mm} 
\frac{3^{dr}}{N_{r}} \frac{(1/3) \cdots (r+1/3)}{r!} < 0.8r^{4/3}.  
\end{equation}

We divide the prime divisors of $D_{r}$ into two sets, according to their size.
We let $D_{r,s}$ denote the contribution to $D_{r}$ from primes less than 
$(3r)^{1/3}$ and let $D_{r,l}$ denote the contribution from the remaining, larger, primes. 

From Lemma~\ref{lem:denom-3}(a), we know that 
\begin{displaymath}
D_{r,s} \leq \prod_{\stackrel{p<(3r)^{1/3}}{p \equiv 2 \bmod 3}}
p^{\lfloor \log(3r)/(2\log(p))+1/2 \rfloor }. 
\end{displaymath}

Now $\lfloor x/2+1/2 \rfloor \leq 3/2 \lfloor x/3 \rfloor +1$, so
\begin{displaymath}
D_{r,s} \leq 
\exp \left\{ \frac{3 \psi \left( \sqrt[3]{3r};3,2 \right)}{2}
             + \theta \left( \sqrt[3]{3r};3,2 \right) \right\}. 
\end{displaymath}

From Theorem~\ref{thm:theta}, and some calculation, we find that
$\theta(x;3,2), \psi(x;3,2) < 0.51x$, so
\begin{equation}
\label{eq:drs2}
D_{r,s} < \exp \left( 1.28 \sqrt[3]{3r} \right).  
\end{equation}

From (\ref{eq:fact}) and (\ref{eq:drs2}), we know that 
\begin{equation}
\label{eq:small}
\frac{D_{r,s} 3^{dr}}{N_{r}} < e^{0.000006r} 
\hspace{3.0mm} \mbox{ and } \hspace{3.0mm} 
\frac{D_{r,s} 3^{dr}}{N_{r}} \frac{(1/3) \cdots (r+1/3)}{r!} 
< e^{0.000006r}, 
\end{equation}
for $r > 200 \cdot 10^{6}$.

We next consider $D_{r,l}$. 

From Lemma~\ref{lem:denom-3}(b), we see that for any positive integer $N$
satisfying $3r/(3N+2) \geq (3r)^{1/3}$, we have 
\begin{displaymath}
D_{r,l} 
\leq \exp \left\{ \sum_{A=0}^{N} \theta(3r/(3A+1);3,2)  
		      - \sum_{A=0}^{N-1} \theta(3r/(3A+2);3,2) \right\}.  
\end{displaymath}

Let $t_{+}(x)$ denote the maximum of $0.5000351$ and $\theta(y;3,2)/y$ 
for all $y \geq x$ and let $t_{-}(x)$ denote the minimum of $0.4999649$ 
and $\theta(y;3,2)/y$ for all $y \geq x$. With the choice $N=200$, 
we can write 
\begin{displaymath}
D_{r,l} 
\leq \exp \left\{ 3r \left( \sum_{A=0}^{200} \frac{t_{+}(600 \cdot 10^{6}/(3A+1))}
						{3A+1} 
				     - \sum_{A=0}^{199} \frac{t_{-}(600 \cdot 10^{6}/(3A+2))}
						  {3A+2} \right) \right\}, 
\end{displaymath}
since $r > 200 \cdot 10^{6}$.

With Theorem~\ref{thm:theta}(b), we calculate the necessary 
values of $t_{+}(x)$ and $t_{-}(x)$ and find that 
\begin{displaymath}
D_{r,l} < e^{0.910993r},  
\end{displaymath}
for $r > 200 \cdot 10^{6}$.

Combining this inequality with (\ref{eq:small}) yields 
\begin{displaymath}
\frac{3^{dr} D_{r}}{N_{r}} < e^{0.911r} 
\hspace{3.0mm} \mbox{ and } \hspace{3.0mm} 
\frac{3^{dr}}{N_{r}} \frac{(1/3) \cdots (r+1/3)}{r!} D_{r} < e^{0.911r}, 
\end{displaymath}
for $r > 200 \cdot 10^{6}$. 

This completes the proof of the lemma. 
\end{proof}

We now need to define our sequence of approximations to $(a/b)^{1/3}$ and find an upper bound on their size.

We start with bounds on the size of the polynomials.

\begin{lemma} 
\label{lem:analytic-bnds}
Let $m, n$ and $r$ be positive integers with $m \leq n/2$ and let $z$ be any real number
satisfying $0 \leq z \leq 1$. Then
\begin{equation}
\label{eq:analytic-ub}
{ \left( 1 + z \right) }^{r} \leq Y_{m,n,r}(z) \leq { \left( 1 + z^{1/2} \right) }^{2r}.
\end{equation}
\end{lemma}

\begin{remark}
The upper bound is best possible as can be seen by considering $z$ near $0$.

For hypergeometric applications, we are particularly interested in $z$ near $1$, where
it appears that the upper bound could be sharpened to
\begin{displaymath}
4^{-r}\frac{(2r)!}{r!}\frac{\Gamma(1-m/n)}{\Gamma(r+1-m/n)}{ \left( 1 + z^{1/2} \right) }^{2r},
\end{displaymath}
although we have been unable to prove this. This is an equality for $z=1$. In the case of $m=1$ and $n=3$,
this extra factor is about $0.8r^{-1/6}$.
\end{remark}

\begin{proof}
We start by proving the upper bound.

We can write 
\begin{eqnarray*}
{ \left( 1 + z^{1/2} \right) }^{2r} 
& = & \sum_{k=0}^{2r} {2r \choose k} z^{k/2} 
\hspace{3.0mm} \mbox{ and } \\
Y_{r}(z) 
& = & \sum_{k=0}^{r} a_{k} z^{k}
= \sum_{k=0}^{r} {r \choose k} 
                 \frac{(r-k+1+m/n) \cdots (r+m/n)}{(1-m/n) \cdots (k-m/n)} z^{k}.  
\end{eqnarray*}

We shall show that 
\begin{displaymath}
{r \choose k} \frac{(r-k+1+m/n) \cdots (r+m/n)}{(1-m/n) \cdots (k-m/n)} z^{k} 
\leq {2r \choose 2k-1} z^{k-1/2} + {2r \choose 2k} z^{k}.
\end{displaymath}

This will prove that $Y_{r}(z) \leq { \left( 1+z^{1/2} \right) }^{2r}$. 

Since $0 \leq z \leq 1$, it suffices to show that 
\begin{equation}
\label{eq:induct}
a_{k} = {r \choose k} \frac{(r-k+1+m/n) \cdots (r+m/n)}{(1-m/n) \cdots (k-m/n)}  
\leq {2r \choose 2k-1} + {2r \choose 2k} = b_{k}. 
\end{equation}

We demonstrate this by induction. 

For $k=0$, (\ref{eq:induct}) holds since $a_{0}$ and $b_{0}$ 
are both equal to $1$. So we can assume that (\ref{eq:induct}) 
holds for some $k$. 

Notice that 
\begin{eqnarray*}
a_{k+1} & = & \frac{(r-k)(r-k+m/n)}{(k+1-m/n)(k+1)}a_{k} \hspace{3.0mm} \mbox{ and } \\ 
b_{k+1} & = & \frac{(r-k)(2r-2k+1)}{(k+1)(2k+1)}b_{k}.
\end{eqnarray*}

Thus
\begin{displaymath}
\frac{a_{k+1}}{b_{k+1}} 
= \frac{(r-k+m/n)}{(r-k+1/2)}
  \frac{(k+1/2)}{(k+1-m/n)}
  \frac{a_{k}}{b_{k}}.  
\end{displaymath}

Since $m \leq n/2$, it is apparent that $(r-k+m/n)/(r-k+1/2) \leq 1$
and that $(k+1/2)/(k+1-m/n) \leq 1$. Since we have assumed that 
$a_{k}/b_{k} \leq 1$, it is also true that $a_{k+1}/b_{k+1} \leq 1$, 
which completes the proof of (\ref{eq:induct}) and hence the upper bound
for $Y_{m,n,r}(z)$.

To establish the lower bound, we again compare coefficients. It is clear that
$a_{0} = {r \choose 0}$ and that $a_{k} \geq {r \choose k}$ for $1 \leq k \leq r$.
Since $0 \leq z \leq 1$, the lower bound holds.
\end{proof}

\begin{lemma} 
\label{lem:pq-ub}
Let $r$ be a positive integer, $a$ and $b$ be positive integers with $b < a$. Put  
\begin{displaymath}
p_{r} = \frac{a^{r} D_{r}}{N_{r}} X_{r}(b/a) 
\hspace{3.0mm} \mbox{ and } \hspace{3.0mm} 
q_{r} = \frac{a^{r} D_{r}}{N_{r}} Y_{r}(b/a). 
\end{displaymath}

Then $p_{r}$ and $q_{r}$ are integers with 
$p_{r}q_{r+1} \neq p_{r+1}q_{r}$ and 
\begin{equation}
\label{eq:qest}
\frac{D_{r}}{N_{r}} (a+b)^{r} \leq q_{r} 
< 1.161 \cdot 10^{39} { \left\{ e^{0.911} 3^{-d} 
		     { \left( a^{1/2} + b^{1/2} \right) }^{2} 
		     \right\} }^{r}. 
\end{equation}
\end{lemma}

\begin{proof}
The first assertion is just a combination of our definitions 
of $D_{r}$ and $N_{r}$ along with an application of 
Lemma~\ref{lem:numer}, while the second one is equation~(16) 
in Lemma~4 of \cite{Baker2}. 

We now prove the upper bound for $q_{r}$.

From Lemma~\ref{lem:analytic-bnds},
\begin{displaymath}
a^{r} Y_{r}(b/a) \leq { \left( a^{1/2} + b^{1/2} \right) }^{2r}.
\end{displaymath}

The upper bound for $q_{r}$ now follows from Lemma~\ref{lem:coeff-ub}(b). 

The lower bound for $q_{r}$ is an immediate consequence of the 
lower bound for $Y_{m,n,r}(z)$ in Lemma~\ref{lem:analytic-bnds}.
\end{proof}

The next lemma contains the relationship that allows the hypergeometic method to provide
good sequences of rational approximations.

\begin{lemma}
\label{lem:relation}
For any positive integers $m,n$ and $r$ with $(m,n)=1$ 
and for any real number $z$ satisfying $0 < z < 1$, 
\begin{equation}
\label{eq:approx}
z^{m/n} X_{m,n,r}(z) - Y_{m,n,r}(z) = (z-1)^{2r+1} R_{m,n,r}(z). 
\end{equation}
\end{lemma}

\begin{proof}
This is (4.2) of \cite{Chud} with $\nu=m/n$.
\end{proof}

We next determine how close these approximations are to $(a/b)^{1/3}$. 

\begin{lemma}
\label{lem:8}
Let $a,b$ and $r$ be positive integers with $b < a$. Then 
\begin{equation}
\label{eq:remest}
\frac{a-b}{200aq_{r}} < \left| q_{r} (a/b)^{1/3} - p_{r} \right| 
< \frac{1.176 \cdot 10^{40} (a-b)}{b} 
  { \left\{ e^{0.911} 3^{-d} { \left( a^{1/2}-b^{1/2} \right) }^{2} 
    \right\} }^{r}. 
\end{equation}  
\end{lemma}

\begin{proof}
Using our definitions of $p_{r}, q_{r}$ and 
$R_{r}(z)$ and the equality expressed in Lemma~\ref{lem:relation}, 
we find that 
\begin{eqnarray*}
\left| q_{r} (a/b)^{1/3} - p_{r} \right| 
& = & \frac{a^{r} D_{r}}{N_{r}} { \left( \frac{a}{b} \right) }^{1/3} 
      { \left( \frac{a-b}{a} \right) }^{2r+1} 
      \frac{(1/3) \cdots (r+1/3)}{(r+1) \cdots (2r+1)} \\
&   & \times {}_{2} F_{1} \left( r+2/3, r+1; 2r+2; (a-b)/a \right). 
\end{eqnarray*}

Since $(a-b)/a$ and the coefficients of this hypergeometric  
function are all positive, we have \\
${}_{2} F_{1} \left( r+2/3, r+1; 2r+2; (a-b)/a \right) > 1$. 
Using the same arguments as in the proof of Lemma~\ref{lem:coeff-ub}, 
we can also show that 
\begin{displaymath}
\frac{(1/3) \cdots (r+1/3)}{(r+1) \cdots (2r+1)} 
> \frac{0.29}{4^{r}r^{1/6}}. 
\end{displaymath}

Combining these inequalities with the lower bound 
for $q_{r}$ in Lemma~\ref{lem:pq-ub}, we obtain 
\begin{equation}
\label{eq:rlbnd}
\left| q_{r} (a/b)^{1/3} - p_{r} \right| 
> { \left( \frac{D_{r}}{N_{r}} \right) }^{2} 
  \frac{0.29(a-b)^{2r}(1+b/a)^{r}}{4^{r}r^{1/6}} 
  \frac{a-b}{aq_{r}}
\end{equation}

Recall that $N_{r}$ is the greatest common factor 
of the numerators of the coefficients of 
$X_{r} \left( 1 - (a-b)z/a \right)$. Since 
$X_{r}(z)$ is a monic polynomial, 
$N_{r} \leq (a-b)^{r}$. The desired lower bound for 
$\left| q_{r} (a/b)^{1/3} - p_{r} \right|$  
now follows from (\ref{eq:rlbnd}) and Lemma~\ref{lem:coeff-ub}(a). 

To obtain the upper bound, we apply Euler's integral 
representation for the hypergeometric function, we have 
\begin{eqnarray*}
\left| q_{r} (a/b)^{1/3} - p_{r} \right| 
& = & \frac{D_{r}a^{r}}{N_{r}} { \left( 1 - \frac{b}{a} \right) }^{2r+1} 
      \frac{(1/3) \cdots (r+1/3)}{r!} { \left( \frac{a}{b} \right) }^{1/3} \\
&   & \times \left| \int_{0}^{1} t^{r}(1-t)^{r} 
				 { \left( 1 - \frac{(a-b)t}{a} \right) 
				 }^{-r-2/3} dt \right|. 
\end{eqnarray*}

Easton (see the proof of his Lemma~8) showed that 
\begin{displaymath}
\left| \int_{0}^{1} t^{r}(1-t)^{r} { \left( 1 - \frac{(a-b)t}{a} \right) 
				   }^{-r-2/3} dt \right|  
\leq (a/b)^{2/3} { \left\{ a { \left( a^{1/2}+b^{1/2} \right) }^{-2} \right\} 
		 }^{r}.  
\end{displaymath}

The lemma now follows from a little algebra and Lemma~\ref{lem:coeff-ub}(b). 
\end{proof}

\section{A Diophantine Lemma}

Finally, we state a lemma which will be used to determine 
an effective irrationality measure from these approximations. 

\begin{lemma}
\label{lem:kappa}
Let $\theta \in \mathbb{R}$. Suppose that there exist $k_{0},l_{0} > 0$ 
and $E,Q > 1$ such that for all $r \in \mathbb{N}$, there are rational 
integers $p_{r}$ and $q_{r}$ with $|q_{r}| < k_{0}Q^{r}$ and 
$|q_{r} \theta - p_{r}| \leq l_{0}E^{-r}$ satisfying 
$p_{r}q_{r+1} \neq p_{r+1}q_{r}$. Then for any rational 
integers $p$ and $q$ with $p/q \neq p_{i}/q_{i}$ for any 
positive integer $i$ and $|q| \geq 1/(2l_{0})$ we have 
\begin{displaymath}
\left| \theta - \frac{p}{q} \right| > \frac{1}{c |q|^{\kappa+1}}, 
\mbox{ where $c=2k_{0}(2l_{0}E)^{\kappa}$ and 
	     $\kappa = \displaystyle \frac{\log Q}{\log E}$.} 
\end{displaymath}             
\end{lemma}

\begin{proof}
In the proof of Lemma~2.8 of \cite{CV}, it is clearly 
noted that this is true. The extra $Q$ which appears 
in the expression for $c$ in the statement of Lemma~2.8 
of \cite{CV} arises only from consideration of the case 
$p/q=p_{i}/q_{i}$ for some positive integer $i$. 
\end{proof}

\section{Proof of Theorem~\ref{thm:main}}

By the lower bound in Lemma~\ref{lem:8}, we need only prove
Theorem~\ref{thm:main} for those rational numbers $p/q \neq 
p_{i}/q_{i}$ for any positive integer $i$. 

All that is required is a simple application of 
Lemma~\ref{lem:kappa} using Lemmas~\ref{lem:pq-ub} 
and \ref{lem:8} to provide the values of 
$k_{0},l_{0},E$ and $Q$. 

From these last two lemmas, we can choose $k_{0} = 1.161 \cdot 10^{39}, 
l_{0} = 1.176 \cdot 10^{40} (a-b)/b$, 
$E = e^{-0.911} 3^{d} { \left( a^{1/2} - b^{1/2} \right) }^{-2}$ 
and $Q = e^{0.911} \cdot 3^{-d} { \left( a^{1/2} + b^{1/2} \right) }^{2}$.

Lemma~\ref{lem:pq-ub} assures us that 
$p_{r}q_{r+1} \neq p_{r+1}q_{r}$. In addition, 
$Q \geq e^{0.911}3^{-1.5}$ $(\sqrt{2}+1)^{2} > 2.78 > 1$ 
and the condition $a > b$ shows that $l_{0} > 0$. 
If $E > 1$ then we can use Lemma~\ref{lem:kappa}. 

The quantity $c$ in Lemma~\ref{lem:kappa} is 
\begin{displaymath}
2.322 \cdot 10^{39} 
{ \left\{ \frac{2.36 \cdot 10^{40} \cdot 3^{d} \left( a^{1/2}+b^{1/2} \right)}
	       {e^{0.911} b \left( a^{1/2}-b^{1/2} \right) }
  \right\} }^{\kappa}. 
\end{displaymath}

Under the assumptions that $a$ and $b$ are positive integers 
with $b < a, E > 1$ and $\kappa < 2$, one can show, by means 
of calculation and arguments from multivariable calculus, that 
$3^{d}e^{-0.911}(\sqrt{a}+\sqrt{b})/(b(\sqrt{a}-\sqrt{b})) < 1.822$,
the maximum occurring for $a=14$ and $b=11$. So we can simplify the 
expression above, bounding it above by 
\begin{displaymath}
2.322 \cdot 10^{39} { \left(4.3 \cdot 10^{40} \right) }^{\kappa}.
\end{displaymath}

By the lower bound in Lemma~\ref{lem:8} for the $p_{i}/q_{i}$'s, we know that the
$c_{1}$ in Theorem~\ref{thm:main} will be a constant times $a$. Furthermore, we
know that, $a \geq 5$ is required in order that $E > 1$ and $\kappa < 2$. So we can
introduce a factor of $a/5$ into our expression for $c$ above, obtaining
\begin{eqnarray*}
2.322 \cdot 10^{39} { \left(4.3 \cdot 10^{40} \right) }^{\kappa}
& < & \frac{10^{40}}{4.3} \frac{a}{5} { \left( 4.3 \cdot 10^{40} \right) }^{\kappa} \\
& < &10^{40}a { \left( \frac{4.3 \cdot 10^{40}}{\sqrt{21.5}} \right) }^{\kappa} \\
& < & 10^{40(\kappa+1)}a,
\end{eqnarray*}
since $\kappa < 2$.

The condition that $E > 1$ (so that $a/2 < b < a$) along with Liouville's
theorem shows that Theorem~\ref{thm:main} is also true if $\kappa \geq 2$. 

By these estimates and Lemma~\ref{lem:kappa} we now know that
Theorem~\ref{thm:main} holds once $|q| \geq 1/(2l_{0}) 
> b/ \left( 2.36 \cdot 10^{40} (a-b) \right)$. There is 
a simple argument we can use to deal with $q$'s of 
smaller absolute value.

If $p/q$ did not satisfy (\ref{eq:result}), then 
$\left| (a/b)^{1/3}-p/q \right| < 1/ \left( 2q^{2} \right)$ 
would certainly hold and $p/q$ would be a convergent in the
continued fraction expansion of $(a/b)^{1/3}$.

Since $b < a$, it follows that $3b^{2/3} < a^{2/3}+(ab)^{1/3}+b^{2/3}$. As a consequence,
$3b^{2/3}\left( a^{1/3} - b^{1/3} \right) < a-b$, or, more conveniently, 
\begin{displaymath}
{ \left( \frac{a}{b} \right) }^{1/3} - 1
= \frac{a^{1/3} - b^{1/3}}{b^{1/3}}  
< \frac{a-b}{3b}. 
\end{displaymath}

So we know that the continued fraction expansion of $(a/b)^{1/3}$ begins
$[1;x, \ldots]$ where $x \geq \lfloor 3b/(a-b) \rfloor$. Therefore $p_{0}=q_{0}=1$
(here $p_{0}/q_{0}$ is the 0-th convergent in the continued fraction expansion
of $(a/b)^{1/3}$), while $q_{1} \geq \lfloor 3b/(a-b) \rfloor$ and it is certainly true
that $q_{1} \geq b/(2.36 \cdot 10^{40}(a-b))$. 

Hence $p/q=1$, in which case $a/b \geq (b+1)/b$ and $E > 1$ imply that
$(a/b)^{1/3}-1 > 1/(4b) \geq 1/(8a)$ and (\ref{eq:result}) holds.

This completes the proof of the Theorem~\ref{thm:main}. 

\section{Proof of Corollary~\ref{cor:main}}

To prove Corollary~\ref{cor:main}, we first need to obtain obtain a
lower bound for $|\sqrt[3]{n} - p/q|$ from the irrationality measure
we have for the appropriate $\sqrt[3]{a/b}$. There are two different 
ways of doing this. 

{\bf (i)} If $\sqrt[3]{a/b}$ is of the form $s \sqrt[3]{n}/t$ then, from  
Theorem~\ref{thm:main}, we obtain 
\begin{displaymath}
\left| \frac{s \sqrt[3]{n}}{t} - \frac{sp}{tq} \right| 
> \frac{1}{c_{1} |tq|^{\kappa+1}}. 
\end{displaymath}
and, as a consequence,  
\begin{displaymath}
\left| \sqrt[3]{n} - \frac{p}{q} \right| 
> \frac{1}{sc_{1} t^{\kappa} |q|^{\kappa+1}}. 
\end{displaymath}

Let us look at the case of $n=2$ to see how we proceed 
here. We have $a=128, b=125, s=4$ and $t=5$. From Theorem~\ref{thm:main},
we have $c_{1}=2 \cdot 10^{97}$ and $\kappa=1.4321$, so 
\begin{displaymath}
\left| \sqrt[3]{2} - \frac{p}{q} \right| 
> \frac{10^{-99}}{|q|^{2.4321}},  
\end{displaymath}
by the above reasoning. 

We wrote a program in Java to calculate the first $500, 000$ partial fractions
and to bound from below the denominators of the first $500, 000$ convergents
in the continued-fraction expansion of $\sqrt[3]{2}$. For this, we used the algorithm
described  by Lang and Trotter \cite{LT} which uses only integer-arithmetic 
and does not require any truncated approximations to $\sqrt[3]{2}$. 

The denominator of the $500, 000$-th convergent is greater than $10^{257, 000}$ 
and it is easy to verify that 
\begin{displaymath}
\frac{10^{-99}}{|q|^{2.4321}} > \frac{0.25}{|q|^{2.4325}} 
\end{displaymath}
for all $q$ whose absolute value is larger than that. 
Thus, it only remains to check that the desired inequality 
is satisfied for all $q$ whose absolute value is at most 
the denominator of the $500, 000$-th convergent. 

Rather than actually checking directly to see if 
\begin{displaymath}
\left| \sqrt[3]{2} - \frac{p}{q} \right| > \frac{0.25}{|q|^{2.4325}} 
\end{displaymath}
held for all these the convergents of $\sqrt[3]{2}$, we simply 
looked at the partial fractions in the following way. 

From the theory of continued-fractions, one can show that 
\begin{displaymath}
\frac{1}{\left( a_{i+1} + 2 \right) q_{i}^{2}} 
< \left| \alpha - \frac{p_{i}}{q_{i}} \right|, 
\end{displaymath}
where $a_{i+1}$ is the $i+1$-st partial fraction 
in the continued-fraction expansion of $\alpha$ 
while $p_{i}/q_{i}$ is the $i$-th convergent. 

As we see in Table~4, the largest partial fraction found for $\sqrt[3]{2}$
was $a_{484708}=4,156,269$. Therefore, the corollary holds for
$|q| > 9 \cdot 10^{13} > ((4156269+2)/4)^{(1/0.4325)}$.
Now a direct check among the smaller convergents completes the proof of
the corollary for $n=2$ (the constant $c_{2}=0.25$ arises here).

We proceeded in the same way for $n=9$, 10, 18, 19, 20, 22, 25, 28, 30, 36, 43, 44, 49, 57,
65, 66, 67, 68, 70, 76, 83 and 84.

{\bf (ii)} The other possibility is that $\sqrt[3]{a/b}$ is of 
the form $s/(t\sqrt[3]{n})$. In this case, we use 
the fact that $|1/x-1/y|= |(x-y)/(xy)|$ and find that 
\begin{displaymath}
\left| \sqrt[3]{n} - \frac{q}{p} \right| 
> \frac{\sqrt[3]{n}}{sc_{1}|p|t^{\kappa}|q|^{\kappa}}. 
\end{displaymath}

We can assume that $|\sqrt[3]{n}-p/q| < 0.5$ and so 
\begin{displaymath}
\left| \sqrt[3]{n} - \frac{q}{p} \right| 
> \frac{\sqrt[3]{n}}
       {sc_{1}t^{\kappa}(\sqrt[3]{n}+1/2)^{\kappa} |p|^{\kappa+1}}. 
\end{displaymath}

We then proceed in the same way as in the previous case. 

It is in this way that we prove the Corollary for $n=3$, 4, 
5, 6, 7, 11, 12, 13, 15, 17, 26, 31, 37, 39, 41, 42, 50, 52,
58, 60, 61, 62, 63, 78, 90, 91, 98 and 100.

\begin{table}[ht]
\begin{center}
\begin{tabular}{||cccrr||} \hline 
$n$   &      $a$             &           $b$           &  $i$      &   $a_{i}$ \\ \hline
$2$   & $2 \cdot 4^{3}$      &         $5^{3}$         & $484,708$ &   $4,156,269$  \\ \hline
$3$   &     $3^{2}$          &         $2^{3}$         &  $13,628$ &     $738,358$  \\ \hline
$4$   & $2 \cdot 4^{3}$      &         $5^{3}$         & $485,529$ &   $8,312,539$  \\ \hline
$5$   & $239645788^{3}$      & $5 \cdot 140145707^{3}$ & $266,405$ &   $3,494,436$ \\ \hline
$6$   & $467^{3}$            &   $6 \cdot 257^{3}$     & $238,114$ &     $466,540$  \\ \hline
$7$   & $44^{3}$             &    $7 \cdot 23^{3}$     & $274,789$ &  $12,013,483$  \\ \hline
$9$   & $9$                  &         $2^{3}$         &  $97,298$ &   $1,063,588$ \\ \hline
$10$  & $5 \cdot 13^{3}$     &     $4 \cdot 14^{3}$    & $371,703$ &   $1,097,381$  \\ \hline
$11$  & $25022^{3}$          &  $11 \cdot 11251^{3}$   & $217,358$ &   $1,352,125$  \\ \hline
$12$  & $9 \cdot 29^{3}$     &    $4 \cdot 38^{3}$     &  $34,767$ &   $1,185,798$  \\ \hline
$13$  & $57^{3}$             &    $13 \cdot 37^{3}$    &  $55,205$ &   $1,406,955$  \\ \hline
$15$  & $5^{2}$              &     $3 \cdot 2^{3}$     & $245,733$ &   $1,571,507$  \\ \hline
$17$  & $18^{3}$             &     $17 \cdot 7^{3}$    & $169,765$ &   $1,536,142$  \\ \hline
$18$  & $9 \cdot 29^{3}$     &    $4 \cdot 38^{3}$     & $300,238$ &   $3,143,844$  \\ \hline
$19$  & $19 \cdot 3^{3}$     &         $8^{3}$         & $138,226$ &     $521,398$  \\ \hline
$20$  & $20 \cdot 7^{3}$     &         $19^{3}$        &  $72,509$ &   $1,840,473$  \\ \hline
$22$  & $11 \cdot 5^{3}$     &     $4 \cdot 7^{3}$     & $232,141$ &     $595,645$  \\ \hline
$25$  & $239645788^{3}$      & $5 \cdot 140145707^{3}$ &  $20,862$ &   $2,449,303$  \\ \hline
$26$  & $3^{3}$             &          $26$           & $252,311$ &   $1,722,109$  \\ \hline
$28$  & $28$                 &         $3^{3}$         & $275,575$ &   $1,654,773$ \\ \hline
$30$  & $10$                 &          $9$            & $228,793$ &     $197,558$  \\ \hline
$31$  & $22^{3}$             &    $31 \cdot 7^{3}$     & $205,544$ &   $1,643,436$  \\ \hline
$36$  & $467^{3}$            &   $6 \cdot 257^{3}$     & $238,549$ &   $2,799,247$  \\ \hline
$37$  & $10^{3}$             &    $37 \cdot 3^{3}$     & $494,731$ &   $6,591,064$  \\ \hline
$39$  & $39^{2} \cdot 2^{3}$ &        $23^{3}$         & $309,275$ &     $483,161$  \\ \hline
$41$  & $100^{3}$            &    $41 \cdot 29^{3}$    & $321,697$ & $417,960,093$ \\ \hline
$42$  & $49$                 &     $6 \cdot 2^{3}$     & $408,968$ &    $409,489$ \\ \hline
$43$  & $43 \cdot 2^{3}$     &         $7^{3}$         & $227,706$ &  $1,359,766$ \\ \hline
$44$  & $44 \cdot 2^{3}$     &         $7^{3}$         & $260,709$ &    $370,994$ \\ \hline
$49$  & $44^{3}$             &    $7 \cdot 23^{3}$     & $273,736$ &  $1,716,211$ \\ \hline
$50$  & $20 \cdot 7^{3}$     &         $19^{3}$        &  $54,577$ &  $2,055,429$ \\ \hline
$52$  & $2 \cdot 2^{3}$      &          $13$           & $379,989$ &  $3,958,641$ \\ \hline
$57$  & $57 \cdot 33^{3}$    &        $127^{3}$        & $110,601$ &    $847,651$ \\ \hline
$58$  & $4 \cdot 2^{3}$      &         $29$            & $172,932$ &    $139,963$ \\ \hline
$60$  & $2 \cdot 2^{3}$      &         $15$            &  $44,247$ &    $461,876$ \\ \hline
$61$  & $4^{3}$              &         $61$            &  $76,517$ &  $3,405,348$ \\ \hline
$62$  & $4 \cdot 2^{3}$      &         $31$            & $400,816$ &    $330,326$ \\ \hline
$63$  & $4^{3}$              &         $63$            & $168,229$ &  $2,664,200$ \\ \hline
$65$  & $65$                 &        $4^{3}$          & $183,363$ & $16,950,688$ \\ \hline
$66$  & $33$                 &    $4 \cdot 2^{3}$      & $179,933$ &    $589,781$ \\ \hline
$67$  & $67$                 &         $4^{3}$         & $419,845$ &    $937,766$ \\ \hline
$68$  & $17$                 &    $2 \cdot 2^{3}$      & $121,095$ &  $1,059,335$ \\ \hline
$70$  & $35$                 &    $4 \cdot 2^{3}$      & $376,116$ &    $582,245$ \\ \hline
$76$  & $19 \cdot 1111^{3}$  &    $2 \cdot 2353^{3}$   & $300,013$ &    $575,574$ \\ \hline
$78$  & $47^{3}$           & $78 \cdot 11^{3}$    & $421,553$ &  $1,145,724$ \\ \hline
$83$  & $83 \cdot 58^{3}$    &         $253^{3}$       & $431,244$ &    $434,543$ \\ \hline
$84$  & $84 \cdot 33856^{3}$ &       $148273^{3}$      & $236,330$ &  $5,018,560$ \\ \hline
$90$  & $3 \cdot 3^{3}$      &     $10 \cdot 2^{3}$    &  $43,615$ &    $314,175$ \\ \hline
$91$  & $9^{3}$              &     $91 \cdot 2^{3}$    & $123,567$ &    $416,579$ \\ \hline
$98$  & $28$                 &         $3^{3}$         & $274,960$ & $23,166,836$ \\ \hline
$100$ & $5 \cdot 13^{3}$     &     $4 \cdot 14^{3}$    & $336,362$ &  $1,383,591$ \\ \hline
\end{tabular}
\caption{Data for the Proof of Corollary~\ref{cor:main}}
\end{center}
\end{table}

\end{document}